\documentclass[11pt]{article}

\usepackage{stdmath}
\usepackage{amsfonts}
\usepackage{amssymb}
\usepackage{amsmath}
\usepackage{fullpage}
\usepackage{graphics}
\usepackage{psfrag}
\usepackage{epsfig}
\usepackage{graphicx}

\title{An Analysis of Primal-Dual Algorithms for Discounted Markov Decision Processes}
\author{Randy Cogill \\ IBM Research Ireland \\ randycog@ie.ibm.com}
\date{\today}

\begin{document}

\maketitle


\begin{abstract}
Several well-known algorithms in the field of combinatorial optimization can be interpreted in terms of the primal-dual method for solving linear programs. For example, Dijkstra's algorithm, the Ford-Fulkerson algorithm, and the Hungarian algorithm can all be viewed as the primal-dual method applied to the linear programming formulations of their respective optimization problems. Roughly speaking, successfully applying the primal-dual method to an optimization problem that can be posed as a linear program relies on the ability to find a simple characterization of the optimal solutions to a related linear program, called the `dual of the restricted primal' (DRP).

This paper is motivated by the following question: What is the algorithm we obtain if we apply the primal-dual method to a linear programming formulation of a discounted cost Markov decision process? We will first show that several widely-used algorithms for Markov decision processes can be interpreted in terms of the primal-dual method, where the value function is updated with suboptimal solutions to the DRP in each iteration. We then provide the optimal solution to the DRP in closed-form, and present the algorithm that results when using this solution to update the value function in each iteration. Unlike the algorithms obtained from suboptimal DRP updates, this algorithm is guaranteed to yield the optimal value function in a finite number of iterations. Finally, we show that the iterations of the primal-dual algorithm can be interpreted as repeated application of the policy iteration algorithm to a special class of Markov decision processes. When considered alongside recent results characterizing the computational complexity of the policy iteration algorithm, this observation could provide new insights into the computational complexity of solving discounted-cost Markov decision processes. 
\end{abstract}

\section{Introduction}
\label{sec:intro}

Markov decision processes (MDPs) are a widely-used model for problems involving sequential decision making under uncertainty. MDPs are used to model the setting where decisions are made in multiple time periods, the decisions made in each period incur some cost, and the decisions made in each period have some influence on the costs that may be incurred in the future. A solution to an MDP is a rule for making decisions in each time period that minimizes some measure of the overall cost incurred over multiple periods. One of the most commonly used cost criteria is the discounted cost incurred over an infinite planning horizon, where costs incurred in the distant future are weighed less heavily than costs incurred in the near future.

Numerous methods exist for computing an optimal decision policy. The most common among these are value iteration \cite{bellman1957}, policy iteration \cite{howard1960}, and linear programming \cite{depenoux1963}. Moreover, there are numerous variants of value iteration and policy iteration that exhibit various computational and performance advantages over the original variants of these algorithms \cite{puterman2009,shlakhter2010}.

Policy iteration is among the most widely-used of these three methods, since it generally requires few iterations to compute an optimal policy in practice. However, despite the widespread use of policy iteration over the past several decades, a clear understanding of its worst-case computational requirements has only been established within the last several years. Collectively, the recent papers \cite{fearnley2010,ye2011,hansen2013,hollanders2012,scherrer2013} have established that policy iteration runs in polynomial time if a certain input parameter (the discount factor) remains fixed across problem instances, but can run in exponential time if the discount factor is allowed to vary for problem instances of increasing size. Hence, policy iteration solves discounted cost MDPs in polynomial time, but not strongly polynomial time \cite{hollanders2012}. Whether there exists a strongly polynomial time algorithm for solving discounted-cost Markov decision processes remains an open question.


Motivated by the search for a strongly polynomial time algorithm for solving discounted-cost MDPs, we introduce and develop a new class of algorithms based on the primal-dual method for solving linear programs. Several well-known algorithms in the field of combinatorial optimization can be interpreted in terms of the primal-dual method. For example, Dijkstra's algorithm, the Ford-Fulkerson algorithm, and the Hungarian algorithm can all be viewed as the primal-dual method applied to the linear programming formulations of their respective optimization problems \cite{papadimitriou1998}. Roughly speaking, successfully applying the primal-dual method to an optimization problem that can be posed as a linear program relies on the ability to find a simple characterization of the optimal solutions to a related linear program, called the `dual of the restricted primal' (DRP).

We will first show that several widely-used algorithms for Markov decision processes can be interpreted in terms of the primal-dual method, where the value function is updated with suboptimal solutions to the DRP in each iteration. We then provide the optimal solution to the DRP in closed-form, and present the algorithm that results when using this solution to update the value function in each iteration. Unlike the algorithms obtained from suboptimal DRP updates, this algorithm is guaranteed to yield the optimal value function in a finite number of iterations. Finally, we show that the iterations of the primal-dual algorithm can be interpreted as repeated application of the policy iteration algorithm to a special class of Markov decision processes. When considered alongside the recent results characterizing the computational complexity of the policy iteration algorithm, this observation could provide new insights into the computational complexity of solving discounted-cost Markov decision processes. 



\section{Overview of the primal-dual method}
\label{sec:PD_overview}

The primal-dual method is a technique for solving linear programs \cite{papadimitriou1998}. Generally speaking, the primal-dual method iteratively updates feasible solutions to a dual linear program, attempting to find a solution that satisfies complementary slackness. Updates to the dual feasible solutions are obtained in each iteration by solving a simpler linear program. 

To be precise, suppose we seek an optimal solution to the following dual form linear program:
\begin{eqnarray}
\label{eqn:dual_LP}
\begin{array}{rl}
\textrm{maximize:} & b^T\lambda \\
\textrm{subject to:} & A^T\lambda \le c.
\end{array}
\end{eqnarray}
In our presentation of the primal-dual method, we will restrict ourselves to the case where $b \ge 0$ and $c \ge 0$. Under these conditions, the linear program (\ref{eqn:dual_LP}) is guaranteed to be feasible and have to have feasible solutions achieving $b^T\lambda \ge 0$. 

A dual feasible solution $\lambda$ is optimal if and only if there exists a primal solution $x$ such that $x$ and $\lambda$ satisfy the complementary slackness conditions
\begin{align*}
(c - A^T\lambda)^Tx &= 0 \\
Ax &= b \\
x &\ge 0.
\end{align*}
In other words, complementary slackness requires that there exist a primal feasible solution with $x_i = 0$ for all constraints such that $(c - A^T\lambda)_i > 0$. 

In each iteration, the primal-dual method checks if a given dual feasible solution satisfies complementary slackness, and generates a dual feasible solution with improved objective value if complementary slackness is not satisfied. For a given dual feasible solution $\lambda$, let $J(\lambda)$ be the set of tight constraints in (\ref{eqn:dual_LP}). That is, 
\[
J(\lambda) = \left\{ j \,\left| \, \sum_{i=1}^m A_{ij}\lambda_i = c_j \right\} \right. ,
\]
where this set could be empty. Throughout this paper we will often express $J(\lambda)$ simply as $J$ for notational compactness, where the dependence on a given dual feasible solution is understood.

Complementary slackness is satisfied if there exists an $x$ satisfying
\begin{align*}
Ax &= b  \\
x &\ge 0 \\
x_j &= 0  \textrm{ for all } j\notin J . \\
\end{align*}
Equivalently, complementary slackness can be verified by solving the linear program
\[
\begin{array}{rll}
\textrm{minimize:} & \mathbf{1}^Ts & \\
\textrm{subject to:} & A\hat{x} + s = b & \\
& \hat{x}_j = 0 & \textrm{for all } j\notin J \\
& \hat{x} \ge 0, s \ge 0.
\end{array}
\]
This linear program is called the \emph{restricted primal} (\textbf{RP}). Note that \textbf{RP} is always feasible since $\hat{x} = 0$, $s = b$ is always a feasible solution. The current dual feasible solution $\lambda$ is optimal if and only if an optimal solution to \textbf{RP} has $\mathbf{1}^Ts = 0$. To simplify notation, we can rewrite \textbf{RP} as
\[
\begin{array}{rl}
\textrm{minimize:} & \mathbf{1}^Ts \\
\textrm{subject to:} & A\hat{x} + s = b \\
& E\hat{x} = \mathbf{0} \\
& \hat{x} \ge 0, s \ge 0,
\end{array}
\]
where $E$ is the $|\overline{J}|\times n$ matrix with $E_{kj} = 1$ if $j$ is the $k$-th element of $\overline{J}$, and $E_{kj} = 0$ otherwise.

Rather than work directly with \textbf{RP}, the primal dual algorithm solves the dual of  \textbf{RP}, which is called \textbf{DRP}:
\[
\begin{array}{rl}
\textrm{maximize:} & b^T\hat{\lambda} \\
\textrm{subject to:} & A^T\hat{\lambda} + E^T\mu \le \mathbf{0} \\
& \hat{\lambda} \le \mathbf{1}.
\end{array}
\]
Since $\mu$ is unrestricted, we can rewrite \textbf{DRP} as
\[
\begin{array}{rll}
\textrm{maximize:} & b^T\hat{\lambda} & \\[2mm]
\textrm{subject to:} & \sum_{i=1}^m A_{ij}\hat{\lambda}_i \le 0 & \textrm{for all } j\in J \\[1mm]
& \hat{\lambda} \le \mathbf{1}. &
\end{array}
\]
The advantage of working directly with \textbf{DRP} is that an optimal $\hat{\lambda}$ can be used to improve the current dual feasible solution when complementary slackness is not satisfied. That is, if complementary slackness is not satisfied, then an optimal solution to \textbf{DRP} has $b^T\hat{\lambda} > 0$. Replacing the current dual feasible solution with $\lambda + \theta\hat{\lambda}$ yields
\[
b^T(\lambda + \theta\hat{\lambda}) > b^T\lambda
\]
for $\theta > 0$. Also, note that
\[
A^T(\lambda + \theta\hat{\lambda}) = A^T\lambda + \theta A^T\hat{\lambda}.
\]
Since 
\[
\sum_{i=1}^m A_{ij}\hat{\lambda}_i \le 0
\]
for all $j$ such that
\[
\sum_{i=1}^m A_{ij}\lambda_i = c_j,
\]
there exists $\theta > 0$ such that $\lambda + \theta\hat{\lambda}$ is dual feasible. In particular, the greatest $\theta$ that maintains dual feasibility is
\[
\theta = \min_{j \in K} \left\{ \frac{ c_j - \sum_{i=1}^m A_{ij}\lambda_i }{ \sum_{i=1}^m A_{ij}\hat{\lambda}_i }  \right\},
\]
where $K = \left\{ k \,\left| \, \sum_{i=1}^m A_{ik}\hat{\lambda}_i > 0 \right. \right\}$. 

\vspace{5mm}

\noindent The primal-dual method can be summarized as follows:
\begin{enumerate}
\item Select an initial dual feasible solution $\lambda$ (such as $\lambda = 0$, when $c \ge 0$).
\item Let 
\[
J = \left\{ j \,\left| \, \sum_{i=1}^m A_{ij}\lambda_i = c_j \right\} \right. .
\]
\item Solve the \textbf{DRP}
\[
\begin{array}{rll}
\textrm{maximize:} & b^T\hat{\lambda} & \\[2mm]
\textrm{subject to:} & \sum_{i=1}^m A_{ij}\hat{\lambda}_i \le 0 & \textrm{for all } j\in J \\[1mm]
& \hat{\lambda} \le \mathbf{1}. &
\end{array}
\]
\item If the optimal achievable value of \textbf{DRP} is $b^T\hat{\lambda} = 0$, then the current dual feasible solution is optimal. Otherwise, update $\lambda := \lambda + \theta\hat{\lambda}$, where
\[
\theta = \min_{j \in K} \left\{ \frac{ c_j - \sum_{i=1}^m A_{ij}\lambda_i }{ \sum_{i=1}^m A_{ij}\hat{\lambda}_i }  \right\}
\]
and $K = \left\{ k \,\left| \, \sum_{i=1}^m A_{ik}\hat{\lambda}_i > 0 \right. \right\}$, then return to Step 2.
\end{enumerate}

We will conclude this section with three important observations regarding the primal-dual method:

\vspace{5mm}

\noindent \textbf{Observation 1:} The primal-dual method solves a linear program by solving a sequence of related linear programs. By doing so, it may seem that the primal-dual method creates more problems than it initially set out to solve. However, for certain problems, \textbf{DRP} is considerably easier to solve than the original linear program. A number of classic combinatorial optimization algorithms with this property are analyzed in \cite{papadimitriou1998}. In Section~\ref{subsec:Optimal_PD}, we will show that the \textbf{DRP} for discounted cost Markov decision processes can be easily solved under certain conditions.

\vspace{5mm}

\noindent \textbf{Observation 2:} We can use the primal-dual method to obtain practical algorithms without necessarily solving \textbf{DRP} to optimality. That is, any feasible solution to \textbf{DRP} that has $b^T\hat{\lambda} > 0$ can be used to obtain a strict improvement in $\lambda$ in each iteration.

\vspace{5mm}

\noindent \textbf{Observation 3:} If \textbf{DRP} is solved to optimality in each iteration, then under reasonable conditions the primal-dual method will converge in finitely many iterations. This comment is made precise in the Lemma below. 

\vspace{1cm}

\noindent \textbf{Lemma 1 \cite{papadimitriou1998}:} The primal-dual method will solve a linear program in finitely many iterations if optimal solutions to \textbf{DRP} are used as updates, and the optimal solutions to \textbf{DRP} are unique in each iteration.

\vspace{1cm}

\noindent \textbf{Proof:} Let $\hat{\lambda}^*(k)$ be the unique optimal solution to \textbf{DRP} and $J(k)$ be the set of tight constraints in iteration $k$. To prove finite convergence, we will show
\[
b^T\hat{\lambda}^*(k+1) < b^T\hat{\lambda}^*(k)
\]
for all iterations $k$ such that $b^T\hat{\lambda}^*(k) > 0$. Since the optimal value of \textbf{DRP} is determined entirely by $J(k)$, this implies that there are no distinct iterations $k \ne k'$ such that $J(k) = J(k')$. Since there are finitely many possible subsets of constraints, this will show that the primal-dual method terminates after finitely many iterations.

If $b^T\hat{\lambda}^*(k) > 0$, then any $j \in J(k)$ such that
\[
\sum_{i=1}^m A_{ij}\hat{\lambda}^*_i(k) < 0 
\]
will not be in $J(k+1)$ in the subsequent iteration. Let $H(k)$ denote the set of constraints
\[
H(k) = \left\{ j\in J(k) \,\left|\, \sum_{i=1}^m A_{ij}\hat{\lambda}^*_i(k) = 0 \right\} \right. .
\]
Since the constraints in $J(k)\setminus H(k)$ are not active, $\hat{\lambda}^*(k)$ is also the unique optimal solution to 
\begin{eqnarray}
\label{eqn:h_LP}
\begin{array}{rll}
\textrm{maximize:} & b^T\hat{\lambda} & \\[2mm]
\textrm{subject to:} & \sum_{i=1}^m A_{ij}\hat{\lambda}_i \le 0 & \textrm{for all } j\in H(k) \\[1mm]
& \hat{\lambda} \le \mathbf{1}. &
\end{array}
\end{eqnarray}
The optimal solution to (\ref{eqn:h_LP}) is unique because otherwise there would exist a convex combination of solutions distinct from $\hat{\lambda}^*(k)$ that was optimal for \textbf{DRP} at iteration $k$.

In addition to constraints departing $J(k)$, at least one $j\notin J(k)$ such that
\[
\sum_{i=1}^m A_{ij}\hat{\lambda}^*_i(k) > 0 
\]
will enter $J(k+1)$ in the subsequent iteration. Note that $H(k)$ is a strict subset of $J(k+1)$. Since 
\begin{eqnarray}
\label{eqn:h_LP_2}
\begin{array}{rll}
\textrm{maximize:} & b^T\hat{\lambda} & \\[2mm]
\textrm{subject to:} & \sum_{i=1}^m A_{ij}\hat{\lambda}_i \le 0 & \textrm{for all } j\in J(k+1) \\[1mm]
& \hat{\lambda} \le \mathbf{1}. &
\end{array}
\end{eqnarray}
is more tightly constrained than (\ref{eqn:h_LP}), 
\[
b^T\hat{\lambda}^*(k+1) \le b^T\hat{\lambda}^*(k).
\]
Moreover, since 
\[
\sum_{i=1}^m A_{ij}\hat{\lambda}^*_i(k) > 0,
\]
for some $j \in J(k+1)$, it cannot be the case that $\hat{\lambda}^*(k+1) = \hat{\lambda}^*(k)$. Finally, since any feasible solution to (\ref{eqn:h_LP_2}) is also feasible for (\ref{eqn:h_LP}), there cannot be a feasible solution to (\ref{eqn:h_LP_2}) that is not equal to $\hat{\lambda}^*(k)$ but has objective value $b^T\hat{\lambda}^*(k)$. Therefore, $b^T\hat{\lambda}^*(k+1) < b^T\hat{\lambda}^*(k)$.
\hfill $\blacksquare$

\vspace{1cm}

\vspace{5mm}

In the next section we will apply the primal-dual method to finite-state, finite-action, discounted cost Markov decision processes. 


\section{Primal-dual for discounted cost MDPs}
\label{sec:PB_MDPs}

In this paper we will examine the algorithms that arise when applying the primal-dual method to a linear program associated with a discounted-cost Markov decision process. The objective of a Markov decision process is to choose actions that influence the evolution of a Markov chain in such a way that the process is directed toward favorable states. Specifically, when action $u$ is chosen, the state of the process evolves according to the transition matrix $P(u)$. We use $c(u)$ to denote the vector of costs incurred in each state when action $u$ is chosen. A static state-feedback policy $\mu$ is a rule that determines the action taken when in each state. We occasionally use the notation $P(\mu)$ and $c(\mu)$ to denote the transition matrix and cost vector induced by policy $\mu$. The discounted cost under policy $\mu$ is
\[
v(\mu) = \sum_{t=0}^\infty \bigl( \gamma P(\mu) \bigr)^t c(\mu),
\]
where $\gamma \in [0,1)$ is a discount factor. For this problem, there always exists an optimal policy $\mu$ such that all components of $v$ achieve their minimum possible value. Our goal is to compute an optimal policy. 

To compute an optimal policy, one can compute a solution $v$ to Bellman's equation
\[
v_i = \min_{u\in\mathcal{U}}\left\{ c_i(u) + \gamma \sum_{j=1}^n P_{ij}(u) v_j \right\} \quad \text{for all } i \in \mathcal{S}.
\]
An optimal policy is then obtained by choosing a minimizing action in each state. One method for solving Bellman's equation is to solve the linear program \cite{depenoux1963}
\[
\begin{array}{rll}
\text{maximize:} & \mathbf{1}^Tv & \\
\text{subject to:} & v \le c(u) + \gamma P(u) v & \text{for all } u \in \mathcal{U} .
\end{array}
\]
We will apply the primal-dual method to this linear program.

The overall structure of the primal-dual method applied to discounted cost MDPs is as follows:
\begin{enumerate}
\item Select an initial dual feasible $v$. For example, we can use $v=0$ when $c(u) \ge 0$ for all $u\in\mathcal{U}$.
\item Given a feasible $v$, determine the set
\[
J = \left\{\, (i,u) \,\left|\, v_i = c_i(u) + \gamma \sum_{j=1}^n P_{ij}(u)v_j \,\right\} \right. .
\]
\item For given $J$, select a $\hat{v}$ such that
\[
\begin{array}{ll}
\mathbf{1}^T\hat{v} > 0, & \\[1mm]
\hat{v}_i \le \gamma \sum_{j=1}^n P_{ij}(u) \hat{v}_j & \text{for all } (i,u) \in J, \\[1mm]
\hat{v} \le \mathbf{1}. &
\end{array}
\] 
If no such $\hat{v}$ exists, then $v$ is optimal.
\item If a feasible $\hat{v}$ is found in Step 3, update $v := v + \theta\hat{v}$, where
\[
\theta = \min_{(i,u) \in K} \left\{ \frac{ c_i(u) + \gamma \sum_{j=1}^n P_{ij}(u)v_j - v_i }{ \hat{v}_i - \gamma \sum_{j=1}^n P_{ij}(u) \hat{v}_j }  \right\}
\]
and $K = \left\{ (i,u) \,|\, \hat{v}_i - \gamma \sum_{j=1}^n P_{ij}(u) \hat{v}_j > 0 \right\}$. Return to Step 2.
\end{enumerate}

In the approach described above, note that we have not required the $\hat{v}$ chosen in Step 3 to solve the \textbf{DRP} to optimality. We will show that several well-known algorithms result from choosing suboptimal $\hat{v}$ satisfying $\mathbf{1}^T\hat{v} > 0$. Then, in Section~\ref{subsec:Optimal_PD} we will present an algorithm that does generate updates by solving \textbf{DRP} to optimality.


\subsection{Primal-dual interpretations of existing algorithms}
\label{subsec:existing}

As our first instance of a primal-dual algorithm, we will use an extremely simple choice of $\hat{v}$ in each iteration. Specifically, in each iteration we will use a $\hat{v}$ such that
\begin{itemize}
\item $\hat{v}_i = 1$ for some $i\in\mathcal{S}$ such that $(i,u)\notin J$ for all $u\in\mathcal{U}$.
\item $\hat{v}_k = 0$ for all $k\ne i$.
\end{itemize}
If there is some $(i,u)\in J$ for all $i\in\mathcal{S}$, then the current dual feasible $v$ is optimal. 

Clearly this choice satisfies $\mathbf{1}^T\hat{v} > 0$. Also, $\hat{v}$ is feasible for \textbf{DRP} since 
\begin{eqnarray*}
\hat{v}_k &=& 0 \\
&\le& \gamma P_{ki}(u) \\
&=& \gamma \sum_{j=1}^n P_{kj}(u) \hat{v}_j.
\end{eqnarray*}
for all $(k,u) \in J$. Finally, in each iteration the update $v + \theta \hat{v}$ is generated using
\begin{eqnarray}
\label{eqn:theta}
\theta = \min_{u} \left\{ \frac{ c_i(u) + \gamma \sum_{j=1}^n P_{ij}(u)v_j - v_i }{ 1 - \gamma P_{ii}(u) }  \right\}.
\end{eqnarray}
This algorithm is equivalent to the well-known variant of value iteration called Gauss-Seidel-Jacobi value iteration \cite{shlakhter2010}. This algorithm and its accelerated performance over ordinary value iteration are typically presented as resulting from a splitting of the transition matrices under each action. Such splittings are discussed in general in \cite{puterman2009}. Below we summarize three related, well-known algorithms that each have a primal-dual interpretation.

\vspace{5mm}

\noindent \textbf{Gauss-Seidel-Jacobi value iteration:} The Gauss-Seidel-Jacobi value iteration algorithm cycles among components of the value function, updating each component as
\[
v_i := \min_{u} \left\{ \frac{ c_i(u) + \gamma \sum_{j\ne i} P_{ij}(u)v_j }{ 1 - \gamma P_{ii}(u) }  \right\}.
\]
This is equivalent to updating the value function as $v := v + \theta \hat{v}$ using the $\theta$ and $\hat{v}$ described above. 

\vspace{5mm}

\noindent \textbf{Gauss-Seidel value iteration:} We obtain the related algorithm known as Gauss-Seidel value iteration if we use the same choice of $\hat{v}$ used in Gauss-Seidel-Jacobi value iteration, but update with a suboptimal choice of $\theta$ given by
\[
\theta = \min_{u} \left\{ c_i(u) + \gamma \sum_{j=1}^n P_{ij}(u)v_j - v_i \right\}.
\]
Specifically, this yields componentwise updates to the value function of the form
\[
v_i := \min_{u} \left\{ c_i(u) + \gamma \sum_{j=1}^n P_{ij}(u)v_j  \right\}.
\]
This algorithm differs from ordinary value iteration in that components are updated sequentially, using recently computed values in each subsequent iteration, rather than updating all components in parallel. 

To show that this choice of $\theta$ is suboptimal, note that
\[
c_i(u) + \gamma \sum_{j=1}^n P_{ij}(u)v_j - v_i \le \frac{ c_i(u) + \gamma \sum_{j=1}^n P_{ij}(u)v_j - v_i }{ 1 - \gamma P_{ii}(u) }.
\]
for all $u$. Therefore, it must be the case that
\[
\min_{u} \left\{ c_i(u) + \gamma \sum_{j=1}^n P_{ij}(u)v_j - v_i \right\} \le \min_{u} \left\{ \frac{ c_i(u) + \gamma \sum_{j=1}^n P_{ij}(u)v_j - v_i }{ 1 - \gamma P_{ii}(u) }  \right\}.
\]
Moreover, the inequality is strict when $P_{ii}(u) > 0$ for all $i$ and $u$.

\vspace{5mm}

\noindent \textbf{Value iteration:} Ordinary value iteration can be interpreted as choosing a suboptimal solution to \textbf{DRP}, then updating $v$ using a suboptimal value of $\theta$. Specifically, ordinary value iteration is equivalent to using primal-dual updates with
\[
\hat{v}_i = \min_{u\in\mathcal{U}}\left\{ c_i(u) + \gamma \sum_{j=1}^n P_{ij}(u) v_j \right\} - v_i \\
\]
and $\theta = 1$. It is easily verified that, provided $v$ is feasible, this choice of $\hat{v}$ is feasible for \textbf{DRP}. Moreover, it is worth noting that the $\theta$ chosen according to Step 4 will always satisfy $\theta \ge 1$. Choosing $\theta$ according to Step 4 rather than $\theta=1$ would yield updates that dominate those of ordinary value iteration, provided that both algorithms are initialized with the same dual feasible $v$.

\vspace{5mm}

Despite the accelerated convergence of Gauss-Seidel-Jacobi value iteration over ordinary value iteration, this algorithm generally does not compute an optimal $v$ in finitely many iterations. In the next section we will examine the convergence behavior of this algorithm on a simple example. For this same example, we then illustrate the finite convergence of the primal-dual algorithm that uses optimal solutions to \textbf{DRP} in each iteration.


\subsection{Examples illustrating convergence}
\label{subsec:convergence}

In general, the variants of value iteration discussed in the previous section do not produce the exact optimal value function in finitely many iterations. Moreover, the closeness to optimality of the solution produced after a fixed number of iterations is sensitive to the choice of discount factor. To show this, consider the simple example below:

\vspace{5mm}
\noindent \textbf{Example 1: Convergence of Gauss-Seidel-Jacobi value iteration}
\vspace{2mm}

Here we consider an example with two states and two actions. Let $P(u)$ and $c(u)$ be the transition matrix and cost vector for action $u$. In this example we will use

\[
P(1) = \bmat{0 & 1 \\ 1 & 0} \qquad P(2) = \bmat{1 & 0 \\ 0 & 1} \qquad c(1) = \bmat{1 \\ 2} \qquad c(2) = \bmat{3 \\ 4}.
\]
Let $v(k)$ be the dual feasible solution produced by Gauss-Seidel-Jacobi value iteration at iteration $k$. The algorithm starts with the dual feasible solution
\[
v(0) = \bmat{0 \\ 0}.
\]
Applying the updates of the algorithm, it easy to show that for odd $k>0$ we have
\[
v_1(k) = 1 + \gamma v_2(k-1),
\]
and for even $k>0$ we have
\[
v_2(k) = 2 + \gamma v_1(k-1).
\]
So, for odd $k>0$ we have
\[
v_1(k) = 1 + \frac{\gamma(2+\gamma)}{1-\gamma^2}(1-\gamma^{k-1}),
\]
and for even $k>0$ we have
\[
v_2(k) = \frac{2+\gamma}{1-\gamma^2}(1-\gamma^{k}).
\]
In the limit this algorithm yields
\[
\lim_{k\rightarrow\infty}v(k) = \bmat{1 \\ 0} + \left(\frac{2+\gamma}{1-\gamma^2}\right)\bmat{\gamma \\ 1 }.
\]
As this example shows, Gauss-Seidel-Jacobi value iteration does not generally converge to an optimal dual solution in finitely many iterations. Furthermore, its rate of convergence is affected by the choice of discount factor $\gamma\in[0,1)$. That is, for even $k > 0$ we have
\[
\|v^* - v(k)\|_\infty = \left(\frac{2+\gamma}{1-\gamma^2}\right)\gamma^k.
\]
For given $k$, $\|v^* - v(k)\|_\infty$ can be made arbitrarily large by choosing $\gamma$ arbitrarily close to $1$. In the next example, we will consider a primal-dual algorithm that will converge in a finite number of iterations, independent of discount factor. 


\vspace{5mm}
\noindent  \textbf{Example 2: Convergence of optimal DRP updates}
\vspace{2mm}

Here we will reconsider the example above, now choosing \emph{optimal} solutions to \textbf{DRP} in each iteration. By doing so, we will obtain the optimal solution to the MDP in two iterations. In Section~\ref{subsec:Optimal_PD} we will present generalization of this algorithm and prove that it always converges in finitely many iterations.

As in the previous example, let $v(k)$ be the dual feasible solution produced by the primal-dual algorithm at iteration $k$. Again, the algorithm will start with the dual feasible solution
\[
v(0) = \bmat{0 \\ 0}.
\]
Since $J$ is empty in the first iteration, the optimal solution to \textbf{DRP} in this iteration is
\[
\hat{v} = \bmat{1 \\ 1}.
\]
For this choice of $\hat{v}$, we have $v(1) = v(0) + \theta\hat{v}$ with
\[
\theta = \frac{1}{1-\gamma},
\]
yielding
\[
v(1) = \frac{1}{1-\gamma}\bmat{1 \\ 1}.
\]
At the start of the next iteration we have 
\begin{eqnarray}
\label{eqn:tight_constraint}
v_1(1) = c_1(1) + \gamma \sum_{j=1}^n P_{1j}(1) v_j(1),
\end{eqnarray}
so $J = \{(1,1)\}$. The optimal solution to \textbf{DRP} in this iteration is now
\[
\hat{v} = \bmat{\gamma \\ 1}.
\]
When producing the update $v(2) = v(1) + \theta\hat{v}$, this choice of $\hat{v}$ maintains tightness of the constraint (\ref{eqn:tight_constraint}). Using 
\[
\theta = \frac{1}{1-\gamma^2}
\]
yields the solution
\begin{eqnarray*}
v(2) &=&  \frac{1}{1-\gamma}\bmat{1 \\ 1} + \frac{1}{1-\gamma^2}\bmat{\gamma \\ 1} \\[3mm]
&=& \bmat{1 \\ 0} + \left(\frac{2+\gamma}{1-\gamma^2}\right)\bmat{\gamma \\ 1},
\end{eqnarray*}
which is in fact the optimal value function for the MDP. Moreover, the optimal value function is computed in two iterations, regardless of the discount factor.

\vspace{3mm}

In the next section we will generalize the approach used in this example, and show that the resulting algorithm always converges in finitely many iterations. In this example, notice that a tight constraint is added for each state and each iteration. If this occurred generally, then we could simply bound the number of iterations by the number of states. That is, if in each iteration at least one constraint for a state became tight while preserving the number of existing tight constraints, the algorithm would terminate in a number of iterations equal to the number of states. As one would expect, the general behavior of this algorithm is not as simple as it appears in this example. In particular, we occasionally encounter the situation where the constraint entering $J$ in an iteration corresponds to a state with a constraint currently in $J$. 


\subsection{Optimal primal-dual updates}
\label{subsec:Optimal_PD}

In this section we present a primal-dual algorithm that uses optimal solutions to \textbf{DRP} in each iteration. To simplify the description of the algorithm, we will introduce some new notation. 

Let $H$ denote a set composed of state-action pairs
\[
H = \{(i_1,u_1),\ldots,(i_{|H|},u_{|H|})\},
\]
where all states appearing in $H$ are distinct. Let $G$ denote the set of states appearing in $H$,
\[
G = \{ i_1,\ldots, i_{|H|} \}.
\]
Let $P_{H,G}$ denote the square matrix with $P_{i_k,i_l}(u_k)$ as its $k,l$ element. Let $\hat{v}_G$ denote the column vector with $\hat{v}_{i_k}$ as its $k$-th element. Let $\overline{G}$ denote the set of states in $\mathcal{S}$ that are not in $G$, 
\[
\overline{G} = \{ j_1,\ldots, j_{|\mathcal{S}|-|G|} \}.
\]
Finally, let $P_{H,\overline{G}}$ denote the possibly non-square matrix with $P_{i_k,j_l}(u_k)$ as its $k,l$ element.

\vspace{5mm}

\noindent The algorithm utilizing optimal \textbf{DRP} solutions is the following:

\vspace{5mm}
\noindent \textbf{Primal-dual algorithm:}

\begin{enumerate}
\item Initialize $v = 0$, $G=\emptyset$, and $H=\emptyset$.
\item If $G = \mathcal{S}$, then the current dual feasible solution $v$ is optimal. Otherwise, proceed to Step 3.
\item Let $\hat{v}_i = 1$ for all $i\in\overline{G}$ and  
\[
\hat{v}_G = \gamma(I - \gamma P_{H,G})^{-1}P_{H,\overline{G}}\mathbf{1},
\]
where $P_{H,G}$ and $P_{H,\overline{G}}$ are constructed as described above.
\item Update $v := v + \theta\hat{v}$, where
\begin{eqnarray}
\label{eqn:b}
\theta = \min_{(i,u) \in K} \left\{ \frac{ c_i(u) + \gamma \sum_{j=1}^n P_{ij}(u)v_j - v_i }{ \hat{v}_i - \gamma \sum_{j=1}^n P_{ij}(u) \hat{v}_j }  \right\}
\end{eqnarray}
and $K = \left\{ (i,u) \,|\, \hat{v}_i - \gamma \sum_{j=1}^n P_{ij}(u) \hat{v}_j > 0 \right\}$.
\item Let $(i,u_1)$ be a state-action pair associated with a constraint achieving the minimum in (\ref{eqn:b}). If $H$ does not contain $(i,u_0)$ for some $u_0\in \mathcal{U}$, add $(i,u_1)$ to $H$ and add $i$ to $G$. Otherwise, remove $(i,u_0)$ from $H$ and add $(i,u_1)$ to $H$.
\item Return to Step 2.
\end{enumerate}

\vspace{5mm}

In the following series of lemmas, we will show that the $\hat{v}$ constructed in Step 3 is an optimal solution to \textbf{DRP}. Lemma 2 proves a property of the algorithm used to show feasibility of $\hat{v}$. Lemma 3 then uses this property to show that $\hat{v}$ is feasible for the \textbf{DRP}. Finally, Lemma 4 proves optimality of $\hat{v}$ by showing that $\hat{v}$ dominates all other feasible solutions to \textbf{DRP}.

\vspace{5mm}

\noindent \textbf{Lemma 2:} Suppose that the minimum $\theta$ in each iteration of the primal-dual algorithm is achieved at a unique state-action pair. Then at each iteration, $J$ contains at most one state-action pair $(i,u_0)$ that is not in $H$. Moreover, the $\hat{v}$ computed in an iteration with $H\subset J$ satisfies
\[
\hat{v}_{i} < \gamma\sum_{ j=1}^n P_{ij}(u_0)\hat{v}
\]
for the state-action pair $(i,u_0)$.

\vspace{1cm}

\noindent \textbf{Proof:} We will prove this claim by induction. At the first iteration, both $H$ and $J$ are empty. 

At the start of a subsequent iteration, first suppose that $H=J$. Assume that the minimum $\theta$ is achieved at a unique state-action pair, say $(i,u_1)$. If there is currently no state-action pair in $H$ containing state $i$, then $(i,u_1)$ is added to both $H$ and $J$ and these sets remain equal in the next iteration. If there is a state-action pair in $H$ containing state $i$, say $(i,u_0)$, then this state-action pair is removed from $H$ and replaced with $(i,u_1)$. Also, $(i,u_1)$ is added to $J$, so $J$ will contain exactly one state-action pair that is not in $H$ in the next iteration.

At the start of a subsequent iteration, now suppose that $J$ contains one state-action pair that is not in $H$, say $(i,u_0)$. If  $(i,u_0)$ is in $J$ but not in $H$, then $(i,u_0)$ was removed from $H$ and replaced with the incoming state-action pair $(i,u_1)$ in the previous iteration. Let $\hat{w}$ be the \textbf{DRP} solution from the previous iteration. Since $(i,u_0)$ was in $H$ and $(i,u_1)$ entered $J$ in the previous iteration, the state-action pairs $(i,u_0)$ and $(i,u_1)$ must satisfy
\begin{eqnarray*}
\hat{w}_{i} &=& \gamma\sum_{ j=1}^n P_{ij}(u_0)\hat{w}_j \\
\hat{w}_{i} &>& \gamma\sum_{ j=1}^n P_{ij}(u_1)\hat{w}_j.
\end{eqnarray*}
Now let $\hat{v}$ be the \textbf{DRP} solution in the current iteration. By Lemma 6 in the appendix,
\[
\hat{v}_{i} < \gamma\sum_{ j=1}^n P_{ij}(u_0)\hat{v}.
\]
So, $(i,u_0)$ will be removed from $J$ in the next iteration. Finally, by an argument identical to the $H = J$ case above, a new state action pair will be added to $J$ and will either be added to $H$ or exchanged for an existing state-action pair in $H$. Therefore, $J$ will contain at most one state-action pair that is not in $H$ at the start of the next iteration. \hfill $\blacksquare$

\vspace{1cm}

\noindent \textbf{Lemma 3:} The $\hat{v}$ constructed in Step 4 of the primal-dual algorithm is a feasible solution to \textbf{DRP}.

\vspace{1cm}

\noindent \textbf{Proof:} The vector $\hat{v}$ is feasible for \textbf{DRP} if and only if $\hat{v} \le \mathbf{1}$ and
\begin{eqnarray}
\label{eqn:drp_feasible}
\hat{v}_i \le \gamma\sum_{ j=1}^n P_{ij}(u)\hat{v}_j
\end{eqnarray}
for all $(i,u)\in J$. 

First we will show that $\hat{v}$ satisfies (\ref{eqn:drp_feasible}). The $\hat{v}$ chosen in the primal-dual algorithm has
\[
\hat{v}_i = \gamma\sum_{ j=1}^n P_{ij}(u)\hat{v}_j
\]
for all $(i,u)\in H$. If $H = J$, then (\ref{eqn:drp_feasible}) is satisfied. If $H$ is a proper subset of $J$, then Lemma 2 shows that there is exactly one state-action pair $(i_0,u_0)$ in $J$ that is not in $H$. Moreover, Lemma 2 shows that
\[
\hat{v}_{i_0} < \gamma\sum_{ j=1}^n P_{{i_0}j}(u_0)\hat{v}_j, 
\]
so (\ref{eqn:drp_feasible}) is satisfied.

To finish the proof, we will show $\hat{v} \le \mathbf{1}$. The $\hat{v}$ chosen in the primal-dual algorithm has $\hat{v}_{\overline{G}} = \mathbf{1}$. To show $\hat{v}_G \le \mathbf{1}$, suppose instead that $\max_{j\in G}\{\hat{v}_j\} > 1$. This, together with inequality (\ref{eqn:drp_feasible}), implies
\begin{eqnarray*}
\hat{v}_i &\le&  \gamma\sum_{ j=1}^n P_{ij}(u)\hat{v}_j \\
&\le& \left(\gamma\sum_{ j = 1}^n P_{ij}(u)\right)\max_{j\in G}\{\hat{v}_j\} \\
&=& \gamma \max_{j\in G}\{\hat{v}_j\} \\
&<& \max_{j\in G}\{\hat{v}_j\}.
\end{eqnarray*}
for all $(i,u) \in H$. However, this is impossible since it implies
\[
\max_{j\in G}\{\hat{v}_j\} < \max_{j\in G}\{\hat{v}_j\}.
\]
\hfill $\blacksquare$

\hspace{1cm}

\noindent \textbf{Lemma 4:} The $\hat{v}$ constructed in Step 4 of the primal-dual algorithm is the unique optimal solution to \textbf{DRP}.

\vspace{1cm}

\noindent \textbf{Proof:} Suppose $\hat{w}$ is an arbitrary feasible solution to \textbf{DRP}. By inequality (\ref{eqn:drp_feasible}), any feasible $\hat{w}$ satisfies
\[
(I - \gamma P_{H,G})\hat{w}_G \le \gamma P_{H,\overline{G}}\hat{w}_{\overline{G}}.
\]
Since $(I - \gamma P_{H,G})^{-1}$ is element-wise nonnegative and $\hat{w}_{\overline{G}} \le \hat{v}_{\overline{G}} = \mathbf{1}$, 
\begin{eqnarray*}
\hat{w}_G &\le& \gamma (I - \gamma P_{H,G})^{-1}P_{H,\overline{G}}\hat{w}_{\overline{G}} \\
&\le& \gamma (I - \gamma P_{H,G})^{-1}P_{H,\overline{G}}\mathbf{1}.
\end{eqnarray*}
Since $\hat{v}_G$ is given by the right-hand side of this inequality, $\hat{w}_G \le \hat{v}_G$ for all feasible $\hat{w}$. Therefore, $\hat{w} \le \hat{v}$ for all feasible $\hat{w}$, implying $b^T\hat{w} < b^T\hat{v}$ for all feasible $\hat{w}\ne\hat{v}$. \hfill $\blacksquare$

Note that the presentation of the algorithm, and subsequent proofs of its properties, assume that a single state-action pair is introduced into $H$ in each iteration. That is, the presentation and analysis appear to disregard the case where multiple constraints become tight simultaneously in a given iteration. It is easy to show that the algorithm still terminates in a finite number of iterations if we simply add one of these state-action pairs, selected arbitrarily, to $H$. This is most easily understood in terms of the primal-dual algorithm's connection to the policy iteration algorithm, which we will elaborate on in the next section.


\section{Preliminary analysis of the primal-dual algorithm}
\label{sec:analysis}

In this section we will provide an initial analysis of the number of iterations required by the primal-dual algorithm. While we do not have a complete characterization of the complexity of the primal-dual algorithm at this time, we will identify a promising direction for further analysis. In particular, the primal-dual algorithm can be interpreted as applying the policy iteration algorithm to a collection of subproblems related to the original MDP. If we can bound the number of iterations required by policy iteration to solve this particular class of subproblems, then we can obtain a bound on the number of iterations required by the primal-dual algorithm. 


We will start by highlighting the connection between the primal-dual algorithm and the policy iteration algorithm. Recall that each iteration of the primal-dual algorithm results in one of two outcomes: Either a new state is added to $G$ or a new state-action pair is added to $H$ for a state already in $G$. Since there may only be $n$ iterations that add a new state to $G$, the difficulty in analyzing the algorithm lies in bounding the number of iterations \emph{between} additions to $G$. It turns out that the iterations of the primal-dual algorithm performed between additions to $G$ are equivalent to iterations of the policy iteration algorithm applied to a particular subproblem.

If a new state is not added to $G$ in a given iteration, then there exists some state $i \in G$ and action $u$ with
\[
\hat{v}_i > \gamma \sum_{j\in G} P_{ij}(u) \hat{v}_j + \gamma \sum_{j\notin G} P_{ij}(u),
\]
and where the constraint in the original MDP for this state-action pair becomes tight. This state-action pair is used to update the policy specified by $H$, and $\hat{v}_{G}$ is recomputed for the new policy in the subsequent iteration. 

Note that the process described above is equivalent to applying the sequential improvement form of the policy iteration algorithm \cite{denardo1982} to a particular subproblem, which is itself a Markov decision process. This subproblem is characterized by a controlled Markov process on the states in $G$. No cost is incurred for state transitions within $G$. Upon transitioning to a state in $\overline{G}$, a cost of $\gamma$ is incurred and the process terminates. Loosely speaking, the aim of this subproblem is to determine the actions to take at states in $G$ as to maximize the amount of time spent within $G$. Due to the finite convergence of policy iteration, the primal-dual algorithm will only perform a finite number of iterations before obtaining a $\hat{v}$ such that
\[
\hat{v}_i \le \gamma \sum_{j\in G} P_{ij}(u) \hat{v}_j + \gamma \sum_{j\notin G} P_{ij}(u),
\]
for all $i\in G$ and $u \in \mathcal{U}$. If this is the case, a new state must be added to $G$ at the end of the current iteration. Note that this provides a worst-case characterization of the number of iterations required before a new state is added to $G$, and the primal-dual algorithm may add a new state to $G$ before policy iteration solves the subproblem on $G$ to completion.


This subproblem is very similar to the well-known \emph{first-passage problem} \cite{derman1970,eaton1962,whittle1982,bertsekas1991}. In the first-passage problem, the aim is to optimize the total cost incurred before transitioning into a terminal set of states. Our subproblem is somewhat simpler than the general first-passage problem since a nonzero cost is only incurred when transitioning to a terminal state, the same cost is incurred for any transition into a terminal state, and cost is discounted.


Let $f_{\text{PI-FP}}(n,m,\gamma,L)$ denote the maximum number of iterations required to solve the first-passage subproblem by sequential improvement policy iteration, where $n$, $m$, $\gamma$, and $L$ are the number of states, number of actions, discount factor, and the number of bits required to specify the transition probabilities, respectively. From the discussion above, the number of iterations required by the primal-dual algorithm can be bounded as
\[
f_{\text{PD}}(n,m,\gamma,L) \le \sum_{k=1}^n f_{\text{PI-FP}}(k,m,\gamma,L),
\]
where $f_{\text{PD}}$ denotes the number of iterations required by the primal-dual algorithm to solve an MDP. That is, the primal dual algorithm requires a number of iterations no greater than the total number of iterations required to solve $n$ first-passage subproblems by policy iteration.


At this point, we can invoke a number of existing results regarding the complexity of policy iteration to gain some further insight. Since all variants of policy iteration provide a strict improvement in the discounted cost from some initial state in every iteration, it is clear that $f_{\text{FP}}(k,m,\gamma,L) \le m^k$. In \cite{melekopoglou1994,littman1995}, an example is provided where sequential improvement policy iteration requires a number of iterations that scales exponentially in the number of states. However, this negative result relies on a particularly poor method for choosing action updates in each iteration. When action updates are selected by choosing the $(i,u)$ with maximum 
\[
v_i - \left( c_i(u) + \gamma \sum_{j=1}^n P_{ij}(u) v_j \right)
\]
in each iteration, a recent analysis of sequential improvement policy iteration \cite{scherrer2013} shows that no more than
\[
n^2(m-1)\left( 1 + \frac{2}{1-\gamma}\log\left(\frac{1}{1-\gamma}\right)\right)
\]
iterations are required to solve any $n$-state, $m$ action discounted-cost MDP. This result is an improvement on a series of recent analyses of \cite{ye2011,hansen2013}. Moreover, the analysis in \cite{hollanders2012}, based on the results for total-cost MDPs of \cite{fearnley2010}, shows that policy iteration can require an exponential number of iterations when the discount factor is allowed to vary with the number of states. 

These recent results for policy iteration provide a promising direction for analyzing the number of iterations required by the primal-dual algorithm. Specifically, we are currently working to resolving the following questions:
\begin{itemize}
\item When applying policy iteration to the first-passage subproblem, the primal-dual algorithm updates the policy specified by $H$ by selecting a state-action pair achieving the minimum in (\ref{eqn:b}). Do the results of \cite{ye2011,hansen2013,scherrer2013} still hold under this selection rule? If so, this would immediately imply that 
\[
f_{\text{PD}}(n,m,\gamma,L) \le \frac{1}{3}(n+1)^3(m-1)\left( 1 + \frac{2}{1-\gamma}\log\left(\frac{1}{1-\gamma}\right)\right).
\]
\item Does sequential improvement policy iteration solve the first-passage subproblem in \emph{strongly} polynomial time? That is, given the special structure of the first-passage subproblem, is it solved in a number of iterations polynomial in $n$ and $m$, independent of $\gamma$? If so, this would imply that $f_{\text{PD}}(n,m,\gamma,L)$ is bounded by a polynomial in $n$ and $m$, independent of $\gamma$. 
\end{itemize}


\section{Conclusions}
\label{sec:conclusions}

In this paper we presented a new algorithm for solving discounted cost Markov decision processes based on the primal-dual method. This algorithm utilizes the optimal solutions to a simple linear program, called the DRP, to update the MDP's value function in each iteration. Several variants of the value iteration algorithm can be interpreted in terms of the primal-dual method, where the value function is updated with suboptimal solutions to the DRP in each iteration. We then presented the algorithm that utilizes optimal solutions to the DRP in each iteration. This algorithm bears a close connection to the policy iteration algorithm, and can be interpreted as repeated application of policy iteration to a special class of first-passage problems. When considered alongside recent results characterizing the computational complexity of the policy iteration algorithm, this observation could provide new insights into the computational complexity of solving discounted-cost Markov decision processes. Specifically, the existence of a strongly polynomial time algorithm for solving discounted-cost Markov decision processes remains an open question. The primal-dual method could provide a fruitful path for further exploration for such an algorithm. 


\bibliography{Primal_Dual_MDP_v2}
\bibliographystyle{plain}

\newpage

\section{Appendix}

The following two lemmas are used in support of the optimality proofs provided in Section~\ref{subsec:Optimal_PD}. Specifically, Lemma 6 provides an intuitively reasonable property with an interpretation related to policy iteration. Specifically, suppose that an existing policy is modified by replacing a single action. Further suppose that action $u_1$ is chosen in place of action $u_0$ at state $i$ because it yields a policy improvement. After updating the value function to reflect this change, immediately switching back to action $u_0$ will not yield a policy improvement. 

\vspace{1cm}

\noindent \textbf{Lemma 5:} Suppose $P\in\mathbb{R}^{n\times n}$ is a sub-stochastic matrix and $e\in\mathbb{R}^{n\times 1}$ is the stochastic vector with $e_k = 1$ for some $k\in\{1,\ldots,n\}$. The unique stochastic vector $\pi\in\mathbb{R}^{n\times 1}$ maximizing
\[
f(\pi) = \sum_{t=0}^\infty \gamma^t (\pi^TP^te)
\]
is $\pi = e$.

\vspace{1cm}

\noindent \textbf{Proof:} Let
\[
h = \sum_{t=0}^\infty \gamma^t (P^te),
\]
so that $f(\pi) = \pi^Th$. The vector $h$ satisfies
\begin{eqnarray}
\label{eqn:h_eqn}
h = e + \gamma P h.
\end{eqnarray}
Since $P$ is sub-stochastic and $\gamma\in[0,1)$,
\[
\max_j\{h_j\} > \sum_{j=1}^n \gamma P_{ij} h_j.
\]
for all $i\in\{1,\ldots,n\}$. Since the equation (\ref{eqn:h_eqn}) states that
\[
h_i = \sum_{j=1}^n \gamma P_{ij} h_j
\]
for all $i \ne k$, this implies $h_k > h_i$ for all $i \ne k$. Therefore, the unique stochastic vector $\pi$ maximizing $\pi^Th$ is $\pi = e$. \hfill $\blacksquare$

\vspace{1cm}

\noindent \textbf{Lemma 6:} Suppose $P\in\mathbb{R}^{n\times n}$ is a sub-stochastic matrix, $c\in\mathbb{R}^{n\times 1}$ is an arbitrary vector, and $v\in\mathbb{R}^{n\times 1}$ is the unique solution to
\[
v = c + \gamma P v.
\]
Moreover, suppose $\pi\in\mathbb{R}^{n\times 1}$ is a sub-stochastic vector, $z\in\mathbb{R}$ is an arbitrary scalar, and
\[
v_k > z + \gamma \pi^Tv 
\]
for some $k\in\{1,\ldots,n\}$. Let $\widehat{P}$ be the matrix obtained by replacing the $k$-th row of $P$ with $\pi^T$, $\widehat{c}$ be the vector obtained by replacing the $k$-th element of $c$ with $z$, and $\widehat{v}$ be the unique solution to
\[
\widehat{v} = \widehat{c} + \gamma \widehat{P} \, \widehat{v}.
\]
The vector $\widehat{v}$ satisfies
\[
\widehat{v}_k < c_k + \gamma \sum_{j=1}^n P_{kj} \widehat{v}_j.
\]

\vspace{1cm}

\noindent \textbf{Proof:} Let $e\in\mathbb{R}^{n\times 1}$ be the stochastic vector with $e_k = 1$. To start, note that 
\[
(e - \gamma \pi)^Tv > z
\]
and
\[
(e - \gamma \pi)^T\widehat{v} = z.
\]
Since $I - \gamma P$ is invertible, the vector $(e - \gamma \pi)^T$ can be expressed as a linear combination of the rows of $I - \gamma P$. That is, there exists some $g\in\mathbb{R}^{n\times 1}$ such that
\[
(I-\gamma P)^Tg = e-\gamma \pi
\]
The inequality
\[
(e-\gamma \pi)v > (e-\gamma \pi)\widehat{v}
\]
is equivalent to
\[
g(I-\gamma P)v > g(I-\gamma P)\widehat{v}.
\]
Since all but the $k$-th elements of $(I-\gamma P)v$ and $(I-\gamma P)\widehat{v}$ coincide, the inequality above implies
\[
g_kc_k > g_k\left(\widehat{v}_k - \gamma \sum_{j=1}^n P_{kj} \widehat{v}_j \right)
\]
To complete the proof we must show that $g_k > 0$.

It is known that
\[
(I - \gamma P)^{-1}  = \sum_{t=0}^\infty (\gamma P)^t
\] 
Therefore,
\begin{eqnarray*}
g_k &=& g^Te \\
&=& \left( \sum_{t=0}^\infty \gamma^t (e-\gamma p)^TP^t \right) e \\
&=& \sum_{t=0}^\infty \gamma^t e^TP^te - \gamma \sum_{t=0}^\infty \gamma^t p^TP^te
\end{eqnarray*}
Finally, from Lemma 5 we have
\[
\sum_{t=0}^\infty \gamma^t e^TP^te \ge  \sum_{t=0}^\infty \gamma^t p^TP^te,
\]
which implies $g_k > 0$.  \hfill $\blacksquare$

\end{document}